\documentclass{article}
\usepackage[utf8]{inputenc}
\usepackage{amsmath,amssymb,dsfont,amsfonts}
\usepackage{mathrsfs}
\usepackage{graphics}
\usepackage{graphicx}
\usepackage{color}

\newcommand{\R}{\mathbb{R}}

\newcommand{\proof}{{\noindent \bf Proof:} }
\newcommand{\eop }{ \hfill $\Box$ }

\newtheorem{theorem}{Theorem}[section]
\newtheorem{proposition}[theorem]{Proposition}
\newtheorem{corollary}[theorem]{Corollary}

\newtheorem{definition}[theorem]{Definition}

\usepackage[bottom=2cm,top=2.5cm,left=2.5cm,right=2.5cm]{geometry}

\title{Geometric aspects of Young Integral:\\ decomposition of flows}
\author{Lourival Lima\footnote{E-mail: l190685@dac.unicamp.br. 
Research supported by CAPES 88882.329064/2019-01.} \ \ \  Paulo Ruffino\footnote{E-mail: ruffino@unicamp.br. 
Research partially supported by CNPq 305212/2019-2, FAPESP 2020/04426-6 and 2015/50122-0.}\ \ \   Pedro Catuogno\footnote{E-mail: pedrojc@unicamp.br. 
Research partially supported by CNPq 302.704/2008-6,
480.271/2009-7 and FAPESP 07/06896-5} }

\date{April 2022}

\begin{document}

\maketitle

\begin{abstract} In this paper we study geometric  aspects of  dynamics  generated by Young differential equations (YDE) driven by $\alpha$-H\"older trajectories with $\alpha \in (1/2, 1)$. We present a number of properties and geometrical constructions on this low regularity context: Young It\^o geometrical formula, horizontal lift in principal fibre bundles, parallel transport, covariant derivative, development and anti-development, among others. Our main application here is a geometrical decomposition of flows generated by YDEs according to diffeomorphisms generated by complementary distributions (integrable or not). The proof of existence of this decomposition is based on an Young Itô-Kunita formula for $\alpha$-H{\"o}lder paths proved by Castrequini and Catuogno (Chaos Solitons Fractals, 2022).


\end{abstract}

\section{Introduction}

In this paper we study geometric  aspects of  dynamics  generated by Young differential equations (YDE) driven by $\alpha$-H\"older trajectories with $\alpha \in (1/2, 1]$. More precisely, given a smooth manifold $M$, we focus on geometrical properties of equations of the type:
\begin{equation}\label{Eq principal}
dx_t = X(x_t)\, dZ_t,
\end{equation}
with initial condition $x_0\in M$ at $t=0$, where $x \rightarrow  X(x) \in  \mathcal{L}(\mathbb{R}^d, T_xM)$ is a smooth assignment of $d$ vector fields on $M$ and $Z \in C^{\alpha}([0,T], \mathbb{R}^d)$ is an $\alpha$-H\"older continuous trajectory in $\R^d$. We say that a path $x:[0,T] \rightarrow M$ is a solution of equation (\ref{Eq principal}) if for all test function $f \in C^{\infty}(M; \R)$ we have that
\begin{eqnarray}
\label{interpretation formula}
f(x_t) = f(x_0)+ \int_0^t Xf(x_s)\, dZ_s,
\end{eqnarray}
where $Xf$ is a short term for $\sum Df(x) X(x)e_i$, with  $e_i$'s the elements of the canonical basis of $\R^d$. The last term of equation (\ref{interpretation formula}) is an integral in the Young sense, see e.g. the classical \cite{Young}, or more recent Hairer and Friz \cite{Friz-Hairer}, Gubinelli et al. \cite{Gubinelli}, Lyons \cite{Lyons}, Castrequini and Russo \cite{CR}, Castrequini and Catuogno \cite{Catuogno}, Cong \cite{Cong}, Ruzmaikina 
\cite{Ruzmaikina}, among many others.

\bigskip

In this scenario of low regularity of  trajectories,  the geometric Young It\^o Formula, Theorem \ref{theorem: Ito formula} opens possibility  to many basic geometric constructions on this dynamics. These topics are explored in the next section, where we prove the existence of horizontal lifts on principal fiber bundles with an affine connection. In particular, considering a Riemannian manifold and its orthonormal bundle, parallel transport and covariant derivatives can be established along $\alpha$-H\"older trajectories. Development and anti-development can also be constructed. 

Motivated by the fact that in many kinds of dynamical systems, in order to obtain local or asymptotic parameters of the dynamics, one performs a befitting decomposition of the associated flow, our main example of application of this low regularity techniques in manifolds concerns a decomposition of the associated Young flow.  In fact, many examples of this kind of decomposition appear in the literature related to distinct geometrical or analytical contexts. We mention few of them: given a system in a semi-simple Lie group, we get a lot of information if we decompose the system into each component of the Iwasawa decomposition (see, e.g. in the stochastic context Malliavin and Malliavin \cite{Malliavin}); given a stochastic flow in a Riemannian manifold, one can write this flow (up to some conditions) as a Markovian process in the group of isometries of the manifold composed with a process in the Lie group of diffeomorphisms which fix the initial condition and has derivatives at this point given by an upper triangular matrix, see Ming Liao \cite{Liao}, \cite{Liao1}. Also, given a flow in an $m$-dimensional manifold with a pair of complementary foliation (i.e. locally the manifold and foliations are diffeomorphic to $\R^k \times \R^{m-k}$) then, locally in time and space, a stochastic  flow can be written as a composition of diffeomorphisms which preserve each of these foliations, see \cite{Melo1}, \cite{Melo2}. We are going to make this last example more precise and explore its potential in the Young integral context.

\bigskip

The decomposition of Young flows is allowed thanks to an It\^o-Ventzel-Kunita type formula in this low regularity context, Theorem \ref{theorem: main decomposition} due to Castrequini and Catuogno \cite{Catuogno}. The framework where we apply this formula is a pair of geometrical distributions (involutive, i.e. which generate a foliation, or not). The main result in Section 3 establishes the local decomposition of an Young flow of diffeomorphisms as one component given by a diffeomorphism generated by  vector fields in one distribution and another component given by diffeomorphism generated by the other distribution. Precise definition are given in Section 3.

\bigskip

The article is organized as follows: in the next section we recall basic properties of the Young integral and prove the relevant geometric results we use latter on. In Section 3 we prove the decomposition of flows given complementary distributions. In Section 4 we present examples. Initially linear systems are treated with a pair of foliations given by affine parallel hyperplanes. We present  conditions to the existence of global decomposition for any time in this context. The last example provides explicit calculations in the case of fiber bundles over homogeneous space $G \rightarrow M=G/H$ where $G$ is a Lie group and $H<G$ is a closed subgroup.


\section{Geometric set up}
\subsection{Young differential equation on manifolds}
We recall that for a general metric space $(M,d)$, a curve $\sigma: [0,T] \rightarrow M$  is  $\alpha$-H\"older continuous, with $\alpha > 0$ if there exists a constant $C>0$, such that 
\begin{eqnarray}
d(\sigma(t),\sigma(s)) \leq C|t-s|^{\alpha},
\end{eqnarray}
for all $s, t \in [0,T]$. This concept extends naturally to a Riemaniann manifold, since it carries the well known induced metric $d(x,y)$ given by 
\begin{equation*}\label{riemannian distance}
d(x,y) = \inf \Big\{\int_0^1||\gamma'(t)||dt; \ \gamma: [0,1] \rightarrow M \mbox{ differentiable such that } \gamma(0)= x, \ \gamma(1)= y \Big\}.
\end{equation*}
See e.g. \cite{Manfredo} among many other classical books. Hence, naturally, $\alpha$-H\"older paths are also well defined in Riemannian manifolds. Most of the classical analytic results on this regularity theory also holds for $\alpha$-H\"older paths in a Riemannian manifold. For instance, composition of a differentiable function with an $\alpha$-H\"older trajectory is also an $\alpha$-H\"older path. Particularly, in a geometrical context, for readers convenience we prove the following

\begin{proposition}
\label{Prop: characterization of alpha Holder}   Let $M$ and $N$ be  Riemannian manifolds, $\dim N\geq 1$.
A path $\sigma:[0,T] \rightarrow M$ is $\alpha$-H{\"o}lder continuous  on $M$ if and only if, for all differentiable map $f: M \rightarrow N$, the path $f(\sigma(t))$ is $\alpha$-H{\"o}lder continuous on $N$. 
\end{proposition}
\proof
There are many interesting ways to prove this result. Here, we use an embedding argument. Initially consider that $N$ is an Euclidean space $ \R^n$ and take $\sigma(t)$ an $\alpha$-H\"older trajectory on $M$.
There exists an isometric embedding $i: M \rightarrow \mathbb{R}^d $ for a  sufficiently large integer $d$ (Nash theorem). 
For sake of notation we write $\sigma_t:= \sigma(t)$.

Since $\| i(x)-i(y)\|_{\R^d} \leq d_M (x,y)$ for all $x,y \in M$, we have the following inequalities: 

\[
 \|i(\sigma_t)-i(\sigma_s)\|_{\R^d} \ \leq \ d_M (\sigma_t,\sigma_s) \  \leq \  C |t-s|^{\alpha},
\]
which implies that $i(\sigma_t)$ is $\alpha$-H\"older in $\R^d$. Now, for any differentiable function
$f:M\rightarrow \R^n$, 
use the fact that it can be extended to a differentiable function $\bar{f}: U \rightarrow \R^n$ defined in a tubular neighbourhood $U$ of $i\circ \sigma([0,T])$ in $\R^N$. Hence, $f(\sigma_t)= \bar{f} (i (\sigma_t))$. Since H\"older regularity is preserved by differentiable functions on Euclidean spaces, $f(\sigma_t)$ is $\alpha$-H\"older continuous in $\R^n$. Mind that, in fact, in the compact set $i(\sigma_t)$ the metrics $d_M$, and $\ell_2$ in $\R^d$ are uniformly equivalents, see Lemma 2.2  \cite{Ledesma}. Hence,  $\ell_2$-norm H\"older regularity in $\R^d$ is equivalent to H\"older regularity on $(M, d_M)$.

\bigskip

For a general Riemannian manifold $N$ and a differentiable map $f: M\rightarrow N$, consider  another isometric embedding $i': N \rightarrow \R ^{d'}$ for an integer $d'$ sufficiently large. Then, the last paragraph shows that $i'\circ f(\sigma)$ is $\alpha$-H\"older in $\R^{d'}$. From Lemma 2.2 \cite{Ledesma} we have that 
there exists a positive constant $C_1$ such that 
\[
d_N(f(\sigma_t),f(\sigma_s) ) \leq C_1 \| i(f(\sigma_t)) - i(f(\sigma_t)) \|_ {\R^{d'}} \leq C_2 |t-s|^{\alpha},
\]
for a positive constant $C_2$, which shows that $f(\sigma_t)$ is $\alpha$-H\"older continuous in $N$.

\bigskip


Conversely, suppose that $f(\sigma_t)\in N$ is $\alpha$-H\"older for all differentiable function $f:M\rightarrow N$. Denote the projections of $i(\sigma_t)\in \R^d$ by  $\sigma_t^j:= p_j \circ i(\sigma_t) $ for each $1\leq j\leq d$. 
Let $\varphi: V \rightarrow W\subset N $ be a local parametrization for $N$, with $V$ an open set in an Euclidean space.
There exist another local parametrization obtained from the previous one, just enlarging the domain by homothety, if necessary, which we call again by $\varphi: \tilde{V} \rightarrow W\subset N $ such that the set  
$\{ (x, 0, \ldots, 0); x= \sigma^j(t) \mbox{ for some } t\in [0,T] \} \subset \tilde{V}$ for all $1\leq j\leq d$. 
Consider the differentiable functions  $f_j:M\rightarrow N $ given by $f_j(x):= \varphi (p_j(i(x)), 0, \ldots, 0))$. Then $f_j(\sigma_t):= \varphi(\sigma^j_ t, 0, \ldots, 0))$ is $\alpha$-H\"older by hypothesis. By metric equivalence in compact sets in the domain of the local parametrization, we have that $\sigma^j_t$ is $\alpha$-H\"older for all $1\leq j\leq d$. We conclude that $\sigma (t) \in M$ is $\alpha$-H\"older continuous on $M$. 

\eop

\bigskip

Before we show conditions for existence and uniqueness of solutions for equation (\ref{Eq principal}), we state the main geometric theorem that is a version of Itô's formula for $\alpha$-H{\"o}lder continuous paths. We start with the definition of the Young integral of a real 1-form:

\begin{definition}[Integration of real 1-forms]
\label{Def: integrals of 1-forms}\emph{ Let $N$ be an $n$-dimensional  differentiable manifold with 
$\bigwedge^1(N)$ the space of real $1$-forms. Consider  $\beta \in \bigwedge^1(N)$ and a chart $(U,(y_1, \ldots, y_n))$ in $N$ such  that}
\begin{eqnarray}
\beta = \sum_{i=1}^n\beta_i \, dy^i.
\end{eqnarray} 
\emph{The integral of $\beta$ along an $\alpha$-H\"older path $x:[0,T] \rightarrow N$ is defined by}
\begin{eqnarray}
\int_0^T \ \beta (x_t) \, dx_t = \sum_{i=1}^n\int_0^T\beta_i dx^i_t,
\end{eqnarray}
\emph{where the above integrals are Riemann-Stieltjes integral of $\beta_i$ with respect to the $i$-th coordinate of the path $x_t$. Among others properties, this integration is independent of the local chart, see e.g,  Abraham, Marsden and Ratiu \cite{Abraham} and Ikeda and  Manabe \cite{Ikeda-Manabe}}.
\end{definition}

The integration of real 1-forms above allows one to integrate many tensor fields in a manifold. In particular, if $F:M \rightarrow \R^d$ is a smooth function, the integration 
\[
\int_0^t DF (x_s) \ dx_s
\]
makes sense, looking at each coordinate of $\R^d$. 
Next Theorem is the basic property of the the $\alpha$-H\"older calculus we are treating in this paper.

\begin{theorem}[Young It\^o Formula]
\label{theorem: Ito formula}
Let $M$ and $N$ be Riemannian manifolds. Consider $x \in \mathcal{C}^{\alpha}([0,T], M ) $ and a smooth function  $F : M  \rightarrow N$. Then 
\begin{equation} \label{Eq: Ito formula}
dF(x_t)= DF(x_s)\  dx_s.
\end{equation}
\end{theorem}

\bigskip

\noindent \textbf{Remark:} We highlight that formula (\ref{Eq: Ito formula}) above means that if $\beta$ is a 1-form in $N$ then 
\begin{equation} \label{Eq: Ito formula on 1-forms}
\int_0^t \beta \ dF(x_s)= \int_0^t (dF(x_s))^* \beta \ dx_s.
\end{equation}
In particular, if $N$ is an Euclidean space:
\begin{equation} \label{Eq: Ito formula on Euclidean space}
F(x_t)=F(x_0) + \int_0^t DF(x_s)dx_s.
\end{equation}

\proof  Initially we prove the result for an Euclidean space $N=\R ^{d}$. We use  again the embedding argument from Nash's theorem: there exists a sufficiently large $p \in \mathbb{N}$ such that $M$ can be isometrically embedded into $\mathbb{R}^{m+p}$. Abusing notation, we have $x \in \mathcal{C}^{\alpha}([0,T], \mathbb{R}^{m+p} )$, $F$ is defined in a tubular neighbourhood of the image of $M$  and  $DF(x) \in L(\mathbb{R}^{m+p}, \mathbb{R}^d )$. By Taylor's formula in Euclidean space,
\[
F(x_{t})-F(x_{s})= DF(x_{s}) \cdot (x_{t}-x_{s}) + R(x_{s}, x_{t}), 
\]
with
\[
R(x_{s}, x_{t})= \int_0^1 (1-u)\, \mathrm{Hess}\,  (F)(x_{s}+ u(x_{t}-x_{s})) (x_{t}-x_{s}, x_{t}-x_{s}) \ du.
\]
Since $F$ is smooth, 
\[
\| R(x_{s}, x_{t}) \|
\leq C \| x_{t} -x_{s} \|^2\leq C^{\prime}|t-s|^{2 \alpha}. \]
Let $\pi =\{s_i\}$ be a partition of $[0,T]$. Then 
\begin{equation}\label{eq1}
F(x_t) -F(x_0)= \sum_i F(x_{s_{i+1}})-F(x_{s_i}) = \sum_i DF(x_{s_i})\cdot (x_{s_{i+1}}-x_{s_i})+ \sum_i R(x_{s_i}, x_{s_{i+1}}).
\end{equation}
We have that 
\[
\sum_i \| R(x_{s_i}, x_{s_{i+1}}) \| \leq C^{\prime} \sum_i |s_{i+1}-s_i|^{2\alpha}\leq C^{\prime}T\sup_{i}|s_{i+1}-s_{i}|^{2 \alpha -1}  .
\]
Thus, since $\alpha>1/2$ we have that 
\begin{equation*}
\lim_{|\pi| \rightarrow 0} \sum_i \| R(x_{s_i}, x_{s_{i+1}}) \| =0. 
\end{equation*}
Take the limit $|\pi|\rightarrow 0$ in  equation (\ref{eq1}) and the definition of Stieltjes (Young) integral to finish the proof in this context.

\bigskip

In general, when $N$ is a Riemannian manifold, consider a local chart $\phi$. The previous calculations hold with $\phi\circ F$ whose integration is independent of the coordinate system.


\eop

Note that multidimensional forms of It\^o formula above, integration by parts etc can be obtained from formula  
(\ref{Eq: Ito formula}) considering the manifold $M$ above as appropriate product spaces. We proceed to prove a theorem of existence and uniqueness of solution for equation (\ref{Eq principal}).  

\begin{theorem} \label{Thm: existencia e unicidade}
Given an initial condition $x_0\in M$, there exists a unique maximal solution of the Young differential equation \emph{(\ref{Eq principal})} such that $x(0)=x_0$. Moreover, there exists a flow of (local) diffeomorphisms associated to the solutions.
\end{theorem}
\proof
A simple way to proof the result for local solutions is based on the existence and uniqueness results in the Euclidean space. In fact, given the initial condition $x_0$, let $(U,\Psi)$ be a chart on $M$ with $x_0\in U$. Let $\tilde{X}:= D \Psi(X(\Psi^{-1}(p)))$ be the induced vector field in the image of $\Psi$. The Young differential equation $dy_t = \tilde{X}(y_t)dZ_t$ has a unique solution local solution $ y_t$ with $y_0=\Psi(x_0)$. See e.g. Lejay \cite{Lejay}, Caruana, Lyons and Thierry \cite{Lyons1}, Li and Lyons \cite{Lyons}, Friz and Hairer \cite{Friz-Hairer} and references therein. Take $x_t = \Psi^{-1}(y_t) \subset U$. We claim that $x_t$ is a solution of equation (\ref{Eq principal}). In fact, consider a test function $f \in C^{\infty}(M)$. By Theorem \ref{theorem: Ito formula}, it follows that
\begin{eqnarray*}
f(\Psi^{-1}(y_s)) &=& f(\Psi^{-1}(y_0)) + \int_0^tD(f\circ \Psi^{-1})(y_s)dy_s \\
&=& f(\Psi^{-1}(y_0))+ \int_0^tDf\Psi_*^{-1}\Psi_*X(y_s)dZ_s \\
&=& f(x_0) + \int_0^tXf(x_s)dZ_s.
\end{eqnarray*}
Moreover, the solution $x_t$ does not depend on the choice of local coordinate. In fact, let $(V,\Phi)$ be another chart on $M$, with $x_t \in U\cap V$ and let $z_t$ be the solution of the Young differential equation $\displaystyle{dz_t = \Phi_*X(z_t)dZ_t}$. Then
\begin{eqnarray*}
dy_t &=& \Psi_*X(y_t)dZ_t \\
&=& \Psi_*\Phi_*^{-1}\Phi_*X(z_t)dZ_t \\
&=& \Psi_*\Phi_*^{-1}dz_t.
\end{eqnarray*}
Hence $\displaystyle{y_t = \Psi \Phi^{-1}(z_t)}$ and therefore $\displaystyle{\Phi^{-1}(z_t) = \Psi^{-1}(y_t) = x_t}$. A maximal solution is obtained in the classical way by extending a local solution up to its explosion time. The existence of local flow of (local) diffeomorphisms is also concluded from the Euclidean case using the same local chart argument. 

\eop

\subsection{Horizontal lifts}

Let $\displaystyle{\{P, M, G, \pi\}}$ be a principal fibre bundle with base $M$, structure group $G$ and total space $P$. In this case $M$ is a smooth, connected and paracompact manifold. The projection $\pi$ is taken as $\pi: P \rightarrow M$.  The group $G$ acts freely on $P$ on the right by the action $\displaystyle{R_g: P \rightarrow P}$ defined by $\displaystyle{R_g(u) = ug}$, for $u \in P$ and $g \in G$.  Let $\mathfrak{g}$ be the Lie algebra of $G$, then an  element $A \in \mathfrak{g}$ generates the exponential $\displaystyle{\{\exp tA, \ t \in \mathbb{R}\}}$, which induces a vector field on $P$ by 
$$
A^*u  = \frac{d }{dt} R_{\exp(tA)}u\  \big|_{t=0},
$$
 If $\Gamma^{\infty}(TP)$ is the section of all smooth vector fields on $P$, then the map $A \rightarrow A^*$, from $\mathfrak{g}$ into $\Gamma^{\infty}(TP)$ is a Lie algebra homomorphism. For more details, see  e.g. Shigekawa \cite{shigekawa}, the classical Kobayashi and Nomizu \cite{kobayashi} among many others. The tangent space $TP$ has a naturally defined subspace called the vertical tangent bundle $VTP$
given by $ VT_uP := \ker d\pi_u $ for all $u\in P$. Note that $A^* u\in VT_uP$ for all $A\in \mathfrak{g}$.

\bigskip

A connection in the principal fibre bundle is an assignment of a horizontal subspace $HT_uP$ of $T_uP$ which is the kernel of a $\mathfrak{g}$-valued $1$-form $\omega$ in $P$ with the following properties: 
\begin{itemize}
\item[(i)] (well-behaved vertically) 
$\omega dR_g = \mathrm{Ad}(g^{-1})\omega$, for all  $g \in G$. Here the linear map $\displaystyle{\mathrm{Ad}(g^{-1}): \mathfrak{g} \rightarrow \mathfrak{g}}$ is the derivative at the identity of the adjoint $\mathrm{Ad}(g^{-1}): G \rightarrow G$ defined by $\displaystyle{\mathrm{Ad}(g^{-1})a = g^{-1}a g}$. 
\item[(ii)] (vertical calibration) $\omega(A^*) = A$, where $A^*$ is a vector field on $VTP$.
\end{itemize}

 Such $1$-form $\omega$ is called a connection form in the principal fibre bundle $\{P, M, G, \pi\}$. Moreover, $\omega$  defines the horizontal tangent bundle $HTP$ given by $HT_uP= \ker \omega_u$. Hence,  for all $u \in P$, the tangent space $T_uP$ splits into $HT_uP \oplus VT_uP$ and $dR_g (HT_uP)= HT_{ug}P$.
 
 \bigskip
 
 Now we have the geometric set up to define the horizontal lift of $\alpha$-H{\"o}lder continuous paths. 

 \begin{definition}
Let $x:[0,T] \rightarrow M$ be an $\alpha$-H{\"o}lder continuous path. Consider $u \in P$, with $\pi(u) = x_0$. The horizontal lift of $x_t$ starting at $u$ is a path $\widetilde{x}: [0,T] \longrightarrow P$ such that:
\begin{itemize}
\item[(i)] $\widetilde{x}_0 = u$.
\item[(ii)] $\pi(\widetilde{x}_t)= x_t$ for all $t \in [0,T]$. 
\item[(iii)] $\displaystyle{\int_0^t \omega \ d \widetilde{x}_s = 0}$ for all $t\in [0,T]$.
\end{itemize}
\end{definition}
Next result shows the existence and uniqueness of the horizontal lift for an $\alpha$-H{\"o}lder continous  paths in a manifold. In the proof we apply the same technique used in Kobayashi and Nomizu \cite{kobayashi} and in Shigekawa \cite{shigekawa} where the existence and uniqueness of horizontal lift were proved in the context of $C^1$ paths and semimartingales respectively.  

\begin{theorem} \label{Thm: Horizontal lift}
Given an $\alpha$-H{\"o}lder continuous path $x: [0,T] \rightarrow M$ and an element $u$ in the fibre $\pi^{-1}(x_0)$, there exists (up to a explosion time) a unique horizontal lift $\widetilde{x}: [0,T] \rightarrow P$ with $\tilde{x}_0=u$.
\end{theorem}
\proof

Consider a local trivialization $\phi: \pi^{-1}(U) \rightarrow U \times G$ with $x_0\in U$ and take the $\alpha$-H{\"o}lder path $\nu_t= \phi^{-1} (x_t,e)$. If the horizontal lift of $x_t$ exists at all, it has to be of the form $\widetilde{x}_t = \nu_t a_t$, where $a_t \in G$ is an appropriate path which makes $\widetilde{x}_t$ horizontal and $\nu_0\, a_0 = u$.  

Let $\Psi: P \times G \rightarrow P$ be the right free action of $G$ on $P$. Then, by Theorem \ref{theorem: Ito formula} we have that 
\[
d\tilde{x}_t = \partial_ 1 \Psi (\nu_t, a_t)\ d \nu_t + \partial_2 \Psi (\nu_t, a_t)\  d a_t.
\]
Hence: 
\begin{eqnarray}
\int_0^t\omega\  d\, \widetilde{x}_t &=&\int_0^t (\partial_ 1 \Psi (\nu_t, a_t))^*\omega  \  d \nu_t + \int_0^t(\partial_2 \Psi (\nu_t, a_t))^*\omega \ da_t.  \nonumber\\
&=& \int_0^tR^*_{a_t}\omega \ d\nu_t + \int_0^t \theta \  da_t,   \label{Eq: connection}
\end{eqnarray}
by the vertical calibration of the connection $\omega$, where $\theta$ is the canonical Cartan 1-form given by $\theta_g(dR_{g}A)= A$ for all $g \in G$ and $A \in \mathfrak{g}$. The lift $\tilde{x}_t$ is horizontal if and only if equation (\ref{Eq: connection}) vanishes for all $t\in [0,T]$, i.e. if and only if $\mathrm{Ad}(a_t^{-1})\omega \ d\nu_t = - \theta \ da_t$. Let $F_1, \ldots, F_n$ be  a basis of the right invariant Lie algebra $\mathfrak{g}$. For all $t\in [0,T]$, there exist $\alpha$-H\'older continuous real functions  $\alpha^1_t, \ldots, \alpha^n_ t$, such that:
\begin{eqnarray}
\int_0^t\omega \ d\nu_s = \sum_{i=1}^nF_i\alpha^i_t.
\end{eqnarray}
Using this notation 
we have that a necessary and sufficient condition such that equation (\ref{Eq: connection}) vanishes is that 

\begin{eqnarray*}
\int_0^t\omega\  d\, \widetilde{x}_t &=&  \sum_{i=1}^n\int_0^t \mathrm{Ad}(a^{-1}_t)F_i \ d\alpha^i +  \int_0^t dR_{a_t^{-1}}\  da_t =0,
\end{eqnarray*}
for all $t\in [0,T]$, i.e., trajectory $a_t$ has to satisfy
\[
da_t = -\sum_{i=1}^n dR_{a_t^{}}\   \mathrm{Ad}(a_t^{-1})F_i\ d\alpha^i_t,
\]
with initial condition $a_0$. There exists a unique solution by Theorem \ref{Thm: existencia e unicidade}, hence there exists a unique horizontal lift $\tilde{x}_t$ up to a explosion. Mind that at the border of the local trivialization, one can extend further the solution applying again the same construction above. The maximal solution covers the whole interval $[0,T]$ (by compactness) if there is no explosion in the fibre.

 Note that for initial element in the fibre $a_0 g$,  the horizontal lift is given by $a_tg$.

 \eop

Besides the dynamics and the principal fiber bundle approach presented so far (which are basic to the next Sections), this low regularity It\^o-Young calculus of Theorem \ref{theorem: Ito formula} allows one to develop further geometrical properties. We mention the following three classical geometric aspects:

\bigskip

\noindent \textbf{A. Parallel Transport and covariant derivative:} 
Given a smooth manifold $M$, consider the frame bundle $BM\rightarrow M$ of basis $u:\R^n \rightarrow T_pM$, with $p\in M$, with the structure group $G=Gl(n, \R )$. Last Theorem applied in this context establishes a parallel transport along $\alpha$-H\"older path $x_t\in M$. In fact, given a horizontal lift $u_t$, the parallel transport of a vector $v\in T_{x_0}M$ is obtained by 
\[
\Big/  \hspace{-1.4mm} \Big/_t \, v = u_t \circ u_0^{-1} (v) \in T_{x(t)}M.
\]
It does not depend on the choice of the horizontal lift. Moreover, if we take the orthonormal frame bundle $OM\rightarrow M$ of basis orthonormal basis given by linear isometries $u:\R^n \rightarrow T_xM$, with $x\in M$, with the structure group $G=O(n, \R )$, the parallel transport is also an isometry. 

\bigskip

Covariant derivative can now be defined along an $\alpha$-H\"older path $x_t\in M$. Given a differentiable vector field $Y$, we have that its covariant derivative along $x(t)$ is given by:
\[
D Y (x_t) =  \Big/  \hspace{-1.4mm} \Big/_t \ \  d \  \Big/  \hspace{-1.4mm} \Big/_t^{-1} \ \  Y(x_t). 
\]
where the differentials are interpreted in the sense of Young (Definition \ref{Def: integrals of 1-forms}).

\bigskip

\noindent \textbf{B. Development and anti-development:} Let $M$ be an $m$-dimensional Riemannian manifold, and consider  an $\alpha$-H\"older continuous path $x:[0,T] \rightarrow \R^m $. Take the horizontal operator $H: OM \times  \R^m \rightarrow HTOM$ where $H(u,v)$ is the horizontal lift of $u(v)\in T_ {\pi(u)}M $ up to $HT_u OM$. The development of $x_t$ on $M$ with initial orthonormal frame $u_0$ is obtained from $u_t$, the solution of the YDE:
\[
 d\, u_t = H(u_t, dx_t),
\]
i.e. $\pi(u_t)$ is the development of $x(t)$ on $M$ (rolling without slipping, with initial ``contact plane'' given by $u_0$). On the other hand, the anti-development of an $\alpha$-H\"older continuous path $x:[0,T] \rightarrow M $ is described using its horizontal lift $\tilde{x}_t$ (Theorem \ref{Thm: Horizontal lift}) with initial condition $\tilde{x}_0$:
\[
y_t = \int_0^t \tilde{x}_s^{-1} \ d\, x_s. 
\]
Note that, as expected, $y_t$ depends on the choice of $\tilde{x}_0$. Compare this approach with the classical Brownian motion approach by Eells and Elworthy \cite{Elworthy}, and  the isotropic L\'evy processes approach in Applebaum and Estrade \cite{Applebaum-Estrade}, among many others.

\bigskip

\noindent \textbf{C. Continuous $\alpha$-H\"older paths in $M$ are solutions of Young differential equations:} As established before, solutions of Young equations driven by $\alpha$-H\"older paths on a manifold are also  $\alpha$-H\"older continuous paths. Reciprocally, every $\alpha$-H\"older continuous paths on $M$ is a solution of a Young differential equation (YDE) driven by an $\alpha$-H\"older functions. In fact, take an embedding $i:M\rightarrow \R^{m+p}$ of $M$ into a sufficiently large dimensional Euclidean space. Let $U$ be a tubular neighbourhood with $\pi: U \rightarrow i(M)$ a projection of $U$ into $i(M)$. Given an $\alpha$-H\"older path  $y_t$ on $M$, let $z_t=i(y_t)$. Then $z_t$ is an $\alpha$-H\"older trajectory in $\R^{m+ p}$.  Consider the YDE in $i(M)$:
\[
dx_t = D \pi (x_t)  \, dz_t.
\]
Then $z_t$ is the solution of this YDE with initial condition $x_0=z_0$: just check that the YDE is the differential version of the identity $z_t =\pi (z_t)$, according to Young It\^o formula of Theorem \ref{theorem: Ito formula}. If the projection $\pi$ is orthogonal, as in Elworthy \cite{Elworthy} then the vector fields are gradients of the embedding. In general, the dynamics of other trajectories starting at $x_0\neq y_0$ depends on the embedding and on the projection. This is an interesting topic to be studied further.



\section{Decomposition of flow generated by Young differential equation}

Let $\mbox{Diff}(M)$ be the infinite dimensional Lie group of smooth diffeomorphisms of a compact connected manifold $M$. The Lie algebra associated to $\mbox{Diff}(M)$ is the infinite dimensional space of smooth vector fields on $M$, see e.g. Neeb \cite{Neeb}, Omori \cite{Omori}, among others. The exponential map $\exp \{tY\} \in \mbox{Diff}(M)$ is the associated flow of diffeomorphisms generated by the smooth vector field $Y$. In this context, given an element $\varphi\in \mbox{Diff}(M)$ the derivative of the right translation is given by $R_{\varphi *}Y= Y(\varphi)$ for any smooth vector $Y$. The derivative of left translation $L_{\varphi *}Y= D\varphi (Y)$, and $\mathrm{Ad}(\varphi)Y = \varphi_* (Y (\varphi^{-1}))$.

\bigskip

 In this Lie group notation, a solution flow $\varphi_t$ of an YDE is written as the solution of a right invariant Young differential equation in the Lie group of diffeomorphisms $\mbox{Diff}(M)$: 
\begin{equation}
    \label{Eq: right inv}
d\varphi_t =  R_{\varphi_{t*}}X\, dZ_t.
\end{equation}
Here we abuse notation in the sense that (using the same notation as in  equation \ref{Eq principal}) one can  write 
\[
X\, dZ_t = \sum_{i=1}^d\  X_i \, dZ_t^i,
\]
where $X_j = X(e_j)$ with $e_j$ the elements of the canonical basis.  Hence, equation (\ref{Eq: right inv}) have to be interpret as 
\begin{eqnarray*}
d\varphi_t = \sum_{i=1}^{d} R_{\varphi_{t*}}X_i\, dZ_t^i.
\end{eqnarray*}


Interesting problems arise when one decomposes a (flow of) diffeomorphism  $\varphi \in \mbox{Diff}(M)$, into composition of convenient prescribed components. This kind of decomposition appears in the literature, for example, in Bismut \cite{Bismut}, Kunita \cite{HK} and many others. In particular, it is also relevant when each component of the decomposition belongs to prescribed subgroups of $\mbox{Diff}(M)$, see e.g Melo et al \cite{Melo2}, Catuogno et al \cite{Catuogno}, Iwasawa and non-linear Iwasawa decomposition \cite{Colonius and Ruffino}, Ming Liao \cite{Liao} among many others.

\bigskip

In this section, we explore the Young calculus to proof the existence of a geometrical decomposition of flows generated by  $\alpha$-H{\"o}lder systems $\varphi_t$ given by equation (\ref{Eq principal}). Suppose that locally $M$ is endowed with a pair of regular differentiable distributions:  i.e., every point $x\in M$ has a neighbourhood $U$ and  differentiable mappings $\Delta^1: U \rightarrow Gr_{k}(M)$ and $\Delta^2: U \rightarrow Gr_{{m-k}}(M)$ respectively, where 
$$
Gr_p (M) = \bigcup_{x \in M}Gr_p(T_x M)
$$ 
is the Grasmannian bundle of $p$-dimensional subspaces over $M$, with $1\leq p\leq m$. We assume that $\Delta^1$ and $\Delta^2$ are complementary in the sense that $\displaystyle{\Delta^1(x) \oplus \Delta^2(x) = T_xM}$, for all $x \in U$.   With this notation we define the subgroup of $\mbox{Diff}(M)$ which is generated by a certain distribution $\Delta$ by:


\begin{equation*}
\mbox{Diff}(\Delta, M) = \mbox{cl} \bigl\{ \exp (t_1X_1) \ldots \exp (t_nX_n), \mbox{ with } X_i\in\Delta, t_i\in\mathbb{R}, \forall n\in\mathbb{N} \bigr\}.
\end{equation*}

 Note that if a distribution $\Delta$ is involutive, then each element of the group $\mbox{Diff}(\Delta, M)$ preserves the leaves of the corresponding foliation.  

\bigskip


In particular, in this Section we focus on the subgroups $\mbox{Diff}(\Delta^1, M)$ and $ \mbox{Diff}(\Delta^2, M)$. The main result of  this paper (Theorem \ref{theorem: main decomposition}) establishes a local decomposition of the solution flow $\varphi_t$ into two components: a curve (solution of an autonomous YDE) in $\mbox{Diff}(\Delta^1, M)$ composed with a  non-autonomous path in $\mbox{Diff}(\Delta^2, M)$. 

\begin{definition}
\label{definition: transversality}
We say that an element $\eta \in \mbox{Diff}( M)$ preserves transversality of  $\Delta^1$ and $\Delta^2$ in a neighbourhood $U \subset M$ if  $\eta_{*} \Delta^2 \left(\eta^{-1}(p)\right)\cap \Delta^1(p) = \{0\}$, for all $p \in U$.
\end{definition}

By continuity, for any pair of complementary distributions, there always exists a neighbourhood of the identity $1d \in \mbox{Diff}(M)$ where all elements in this neighbourhood preserve transversality. Moreover, if the distribution $\Delta^1$ is involutive then all elements in $\mbox{Diff}(\Delta^1, M)$ preserves tranversality of $\Delta^1$ and $\Delta^2$: in fact, the derivative $\eta_*$ above is a linear isomorphism which sends tangent spaces of the associated foliation to tangent spaces in the same leaf. In the sequence, we state an extended scope of the It\^o-Kunita formula (see \cite{HK} ) in the geometrical Young calculus.

\begin{theorem}[Young Itô-Kunita formula]\label{Theorem: Ito-Ventzel}
Let $X$, $Y$ $\in C^2(M,\mathcal{L}(\mathbb{R}^d, TM))$ and $Z \in C^{\alpha}([0,T], \mathbb{R}^d])$ and suppose that $\eta_t$ and $\psi_t$ are solutions maps associated to the Young differential equations $d\eta_t = X(\eta_t)dZ_t$ and $d\psi_t=Y(\psi_t)dZ_t$ respectively. Then, $\varphi_t = \eta_t \circ \psi_t$ is the solution map associated with the Young differential equation
\begin{equation}\label{IK}
d\varphi_t = X(\varphi_t)dZ_t + \mathrm{Ad}(\eta_{t}) Y( \varphi_t)dZ_t.
\end{equation}
\end{theorem}

For a proof in this low regularity context, see Castrequini and Catuogno \cite[Thm. 4.1]{RP}. Next Corollary shows that the inverse of the solution flow of an YDE is also $\alpha$-H\"older continuous. 
\begin{corollary}
\label{Corollary: inverse formula}
If $\eta_t$ is the solution flow the Young differential equation on $M$
\begin{eqnarray} \label{Eq: fluxo inverso 1}
dx_t = X(x_t)\, dZ_t,
\end{eqnarray}
then, the inverse map $\eta^{-1}_t$ is the solution of the Young differential equation on $M$
\begin{eqnarray} \label{Eq: fluxo inverso 2}
dz_t = -D\eta^{-1}_t(z_t)X(\eta_t(z_t))dZ_t.
\end{eqnarray}
\end{corollary}

\proof In fact, just apply expressions (\ref{Eq: fluxo inverso 1}) and (\ref{Eq: fluxo inverso 2}) into equation  (\ref{IK}). 

\eop

For a constructive proof of last Corollary see \cite[Thm. 4.2]{RP}. Next Theorem states the main result of this section:

\begin{theorem}[Decomposition of flows of YDE]
\label{theorem: main decomposition}
Up to a life time $\tau \in [0,T]$, the solution flow $\varphi_t$ can be 
locally decomposed as
\[
 \varphi_t = \eta_t \circ \psi_t,
\]
where $\eta_t$ is solution of an (autonomous) Young differential equation in $\emph{Diff}(\Delta^1, M)$ and $\psi_t$ is
a path in $\emph{Diff}(\Delta^2, M)$.
\end{theorem}

\proof


Given $p \in M$, take $\eta \in \mbox{Diff}(\Delta^1,M)$ sufficiently close to the identity such that it preserves tranversality, i.e. $\mathrm{Ad}(\eta_t)\Delta^2$ and $\Delta^1$ are complementary.  The tangent vector(s) $X(p)$ can be decomposed uniquely as 
\begin{eqnarray} \label{Eq: h and V}
X(p) = h(p)+V(\eta_t,p),
\end{eqnarray}
where $h(p) \in \Delta^1(p)$ and $V(\eta_t,p) \in \mathrm{Ad}(\eta_t)\Delta^2(p)$, for all $p \in M$. We take the first component  $\eta_t$ as the solution map of the following Young differential equation in $\mathrm{Diff}(\Delta^1, M)$:
\begin{eqnarray}
\label{equation: horizontal components}
d\eta_t = R_{\eta_{t*}}\, h\, dZ_t,
\end{eqnarray}
with initial condition $\eta_0=1d$, the identity.
Even though the equation above is described in terms of a right translation, it is not a right invariant equation since  $h$ in general depends on $\eta_t$. We obtain the second component of decomposition of $\varphi_t$ using that  $\psi_t = \eta^{-1}_t \circ \varphi_t$.  Applying Corollary \ref{Corollary: inverse formula}, it follows that:
\begin{eqnarray*}
d\eta^{-1}_t = -L_{\eta^{-1}_{t*}}\, h\, dZ_t,
\end{eqnarray*}
where $L_{\eta^{-1}_{t*}}$ is the derivative of the left translation at the identity by $\eta^{-1}_t$. Finally, we find a equation for $\psi_t$ by applying Theorem \ref{Theorem: Ito-Ventzel}:
\begin{eqnarray} \label{Eq: vertical component}
d\psi_t &=& (\eta^{-1}_t h \, \eta_t \, \psi_t - \eta^{-1}_t X\, \eta_t \, \psi_t)\, dZ_t \nonumber\\
&& \nonumber \\
&=& \mathrm{Ad}(\eta^{-1}_t)V(\eta_t)\ dZ_t. 
\end{eqnarray}
Note that $V(\eta, p)$ does not necessarily belong to $\Delta^2$. Still, $d\psi_t \in \Delta^2$ since $d\psi_t \in \mathrm{Ad}(\eta^{-1})\mathrm{Ad}(\eta)\Delta^2 = \Delta^2$. Then $\psi_t$ is the $\Delta^2$-component of $\varphi_t$.

\eop

\begin{corollary}
If the distributions $\Delta^1$ and $\Delta^2$  are integrable, then the decomposition of Theorem \ref{theorem: main decomposition} is unique. 
\end{corollary}
\proof
In fact, in this case $\mathrm{Diff}(\Delta^1, M) \cap \mathrm{Diff}(\Delta^2, M)= \{1d\}$.

\eop

\section{Examples of decomposition of flows of YDE}
In this section we consider the same geometric structure used in section 2, i.e, a principal fibre bundle $\{P, M, G, \pi\}$, with base $M$, structure group $G$ and total space $P$.  Our main goal is to apply the decomposition which was proposed in theorem \ref{theorem: main decomposition} in general fibre bundles.  Important notions such as $1$-forms connection on fibre bundles, horizontal and vertical tangent bundles and others were discussed briefly on section $2$, more details can be found for example in  Kobayashi and Nomizu \cite{kobayashi}.

\subsection{Linear systems}

Consider an Euclidean space $\R^n$, with a pair of complementary foliations given by the trivial Cartesian product $\R^k \times \R^{\ell}$, with 
$k+ \ell = n$. More precisely, the horizontal foliation $\mathcal{F}_H$ is given by parallel leaves generated by  affine translations  $x+ (\R^k \times \{0\}),$ with  $x\in \R^n$. Analogously, the vertical foliation $\mathcal{F}_V$ is given by parallel vertical  leaves $x+ (\{ 0\} \times \R^{\ell})$, for all $x\in \R^n$. We consider the linear Young differential equation:
\begin{equation} \label{Eq: Linear}
  dx_t = A\, x_t \  dZ_t,
\end{equation}
with $x_0\in \R^n $ and $Z_t$ an $\alpha$-H\"older continuous trajectory in the real line. The Young calculus presented in the previous section shows that the  fundamental linear solution flow of (\ref{Eq: Linear}) is the exponential 
\begin{equation} \label{Eq: flow of Linear}
  F_t = \exp{ \{ A (Z_t-Z_0) \} }.
 \end{equation}

\[
A = \left(  \begin{array}{ll}
\Big( A_1 \Big)_{k\times k}     & \Big( A_2 \Big) \\
    \Big( A_3 \Big) & \Big( A_4 \Big)_{\ell \times \ell}
\end{array}
\right)
\]

The decomposition we are interested here is 
\[
F_t = \eta_t \circ \psi_t
\]
such that $\eta_t \in Dif^H$ and $\psi_t \in Dif^V$. In general $\eta_t$ and $\psi_t$ does not have to be linear, even in quite symmetric situations. For example, if the pair of foliations  in $\R^n \setminus \{0\}$ are given by radial and spherical coordinates, the components of the decomposition are not necessarily linear: in fact, the linear radial diffeomorphisms is reduced to a one dimensional group of uniform contractions and expansions $\lambda 1d$, with $\lambda>0$, which, obviously, is not big enough to perform the decomposition. For the Cartesian pair of foliation $\R^k \times \R^{\ell}$ considered in this section, we do have that $\eta_t$ and $\psi_t$ are linear. In fact, in coordinates, write 
\[
F_t = \left(  \begin{array}{ll}
\Big( F_1 (t) \Big)_{k\times k}     & \Big( F_2(t)\Big)_{k \times \ell} \\
    & \\
    \Big( F_3 (t)\Big)_{\ell \times k} & \Big( F_4(t)\Big)_{\ell \times \ell}
\end{array}
\right).
\]
Since $\eta_t$ does not change the last  $\ell$ coordinates the diffeomorphisms $\psi$ must satisfies
\[
\psi_t = \left(  \begin{array}{ll}
   \Big( 1d \Big)_{k\times k}     &  0  \\
   & \\
     F_3(t) & F_4 (t)
\end{array}
\right).
\]
Hence diffeomorphisms  $\psi_t$ and $\eta_t$, when exist, are global and linear.

\bigskip

\noindent \textbf{A simple example:} A system which illustrates not only these formulae, but also the lifetime of the decomposition is the pure rotation in $R^2$ given by 
\[
dx_t = \left( \begin{array}{cc} 0 & -1 \\ 1 & 0 \end{array} \right) x_t\, dZ_t,
\]
whose decomposition of flow can be easily calculated as:
\begin{eqnarray}
\label{Ex: explosion}
\hspace{-1cm}\left( \begin{array}{cc} \cos Z_t & -\sin Z_t \\ \sin Z_t & \cos Z_t \end{array} \right) &=& \left( \begin{array}{cc} \sec Z_t & -\tan Z_t \\ 0 & 1 \end{array} \right)\left( \begin{array}{cc} 1 & 0 \\ \sin Z_t & \cos Z_t \end{array} \right). 
\end{eqnarray}
Note that if 
$Z_{t} \in \{\frac{\pi}{2}+k\pi, k \in \mathbb{Z}\}$, then the decomposition (\ref{Ex: explosion}) no longer exists at the corresponding time t, i.e. we have explosion of the solutions of equations (\ref{equation: horizontal components}) or (\ref{Eq: vertical component}).

\eop

\bigskip

 Back to the general linear case, the components of the decomposition in fact lie in the Lie group:

\[
\psi_t \in G_V= \left\{  g\in Gl(n, \R); g= \left(\begin{array}{ll}
\Big( 1d \Big)_{k\times k}     &  0  \\
     g_3  & \Big( g_4 \Big)_{\ell \times \ell}
\end{array}
\right)  \right\}
\]
whose Lie algebra is given by the vector space generated by 
\[
\left(\begin{array}{cc}
\Big( 0 \Big)_{k\times k}     & 0 \\
& \\
     \Big(  * \Big)  & \Big( * \Big)_{\ell \times \ell}
\end{array}
\right),
\]
where $(*)$ means nonzero matrices of the appropriate dimension. Analogously for the horizontal component: 
\[
\eta_t \in G_H= \left\{  g\in Gl(n, \R); g= \left(\begin{array}{ll}
\Big( g_1 \Big)_{k\times k}     &  g_2  \\
     0  & \Big( 1d \Big)_{\ell \times \ell}
\end{array}
\right)  \right\}
\]
whose Lie algebra is given by the vector space generated by 
\[
\left(\begin{array}{cc}
\Big( * \Big)_{k\times k}     & \Big(  * \Big)_{k \times \ell} \\
& \\
     0  & \Big( 0 \Big)_{\ell \times \ell}
\end{array}
\right),
\]

Using the properties of the Young integral, we find the differential equations for the constituents submatrices $g_1, g_2$ and  $g_3, g_4$ of $\eta_t$ and $\psi_t$ respectively. Let $\pi_2: R^k\times \R ^{\ell}\rightarrow \R ^{\ell}$ be the projection on the second subspace. From formula (\ref{Eq: h and V}) we have that
\[
V(\eta, \cdot) = \eta \circ \pi_2 \circ A (\cdot).
 \]
In fact, it is enough to check that $V(\eta, \cdot)$ is in the image of the vertical component by $\eta$ and that  $\pi_2 V(\eta, \cdot) = \pi_2 A (\cdot)$. From this formula,   equations (\ref{equation: horizontal components}) and (\ref{Eq: vertical component}) we find the autonomous equation:
\[
d\eta_t = (1d - \eta_t \circ \pi_2)A \, \eta_t \ dZ_t,
\]
and the well expected nonautonomous vertical diffeomorphisms:
\[
d\psi_t =   \pi_2\,  A\,  \eta_t \circ \psi_t \ dZ_t.
\]
Rewriting each constituent submatrices we find:
\begin{eqnarray}
\begin{array}{ll}
d g_1(t) = \big[  A_1 \ g_1(t) - g_2(t)\ A_3 \ g_1(t) \big] \, dZ_t,    \ \ \ \ \ \ \ \  &   d g_2(t) = \big[ A_1 g_2(t)+ A_2 - g_2(t)\, A_4 - g_2\, A_3\, g_2(t) \big] \, dZ_t,  \\
  & \\
d g_3(t) = \big[ A_3\,g_1 + A_3\, g_2\, g_3 + A_4 g_3 \big] \, dZ_t,     &   d g_4(t) = \big[ A_3\, g_2\, g_ 4 + A_4\, g_4 \big] \, dZ_t. \label{Eq: constituent g_3 g_4}
\end{array} 
\label{Eq: constituent}
\end{eqnarray}
Explosion in the solutions of the equations of $g_1$ and $g_2$ can appear if $A_3$ is not zero (see example of equation (\ref{Ex: explosion}), where $A_3=[1]$). Otherwise, if $A_3=0$ then there exists the decomposition for all time $t \geq 0$. Using this feature, and the Jordan canonical form we can extend the scope of the decomposition in the next Proposition. Before that let us fix a notation. Given two complementary subspaces $E_1 \oplus E_2 = \R^n$, let us denote by $\mathcal{F}(E_1)$  and $\mathcal{F}(E_2)$ the corresponding pair of complementary parallel foliations in $\R^n$.

\begin{proposition}
Consider a Young linear system in $\R^n$ 
\begin{equation}\label{Eq: linear Young}
dx_t = A\, x_t \ dZ_t.
\end{equation}
If dimension $n>2$, then there exist a pair of parallel foliations $\mathcal{F}(E_1)$, $\mathcal{F}(E_2)$ generated by complementary subspaces $E_1$ and $E_2$ such that the decomposition of the flow of equation (\ref{Eq: linear Young}) exists for all time $t\in [0,T]$, i.e. there is no explosion time of the decomposition. Dimension of $E_1$ can be chosen as a number of the form $(a+2b)$ where $a=0, 1, \ldots, r= \# \{ \mbox{real eigenvalues with multiplicities} \}$, and $b=0, 1, \ldots, (n-r)/2$.
\end{proposition}
\proof
Let $A = P\, J P^{-1}$ be the canonical real Jordan form of $A$, with the choice of bases $P$ such that the nilpotent component has, if  necessary,  $1$'s and identities $I_2$'s above the diagonal. The change of coordinates $y=P\, x$ establishes the conjugate Young system:
\begin{equation*}
dy_t = J\, y_t \ dZ_t.
\end{equation*}
If $n>2$, it is possible to write 
\[
J = \left(  \begin{array}{ll}
\Big( J_1 \Big)_{k\times k}     & \Big( J_2 \Big) \\
    \Big( J_3 \Big) & \Big( J_4 \Big)_{\ell \times \ell}
\end{array}
\right)
\]
with $k=a + 2b$ and its complementary $\ell=n-k$, such that the submatrix $(J_3)_{\ell \times k} = 0$. The number $a$ represents the number of real eingenvalues in the block $J_3$ and $b$ represents the number of pairs of conjugate nonreal eigenvalues in this block. Hence, equations (\ref{Eq: constituent}) guarantee the there is no explosion in the decomposition of $y_t$. By conjugacy, there is also no explosion in the decomposition of the linear fundamental solution $F_t$ of (\ref{Eq: linear Young}) along the foliations generated by $E_1 = P \, (\R^k\times \{0\} )$ and  $E_2 = P \, (\{0\} \times \R^l)$. This proves the Proposition.

\eop

Using the notation in the proof of last Proposition, the decomposition of $F_t= \eta_t \circ \psi_t$ above are such that $\eta_t$ lies in the group $P\, G_H \, P^{-1}$ and $\psi_t$ lies in $P\, G_V\,  P^{-1}$.

\subsection{Principal fibre bundles over homogeneous spaces}

Let $G$ be a connected Lie group with a closed subgroup $H$ and denote by $\mathfrak{g}$ and $\mathfrak{h}$ their Lie algebras of right invariant vector fields, respectively.  The group $G$ acts on $H$ by  left translation $gH$, for all $g \in G$ and the orbits generate the homogeneous space $M:=G/H$, see e.g. \cite{kobayashi}.  We have a principal fibre bundle given by the canonical projection $\pi: G \rightarrow M$. Given an element $A\in \mathfrak{g}$ consider the right invariant YDE:

\begin{eqnarray}
\label{equation: principal fibre bundle}
d\, g_t =  A  g_t\ dZ_t.
\end{eqnarray}

As it has done in Section 2.2, here, we consider a connection $\omega$ in the principal fibre bundle $\pi: G \rightarrow M$. In this example we construct our decomposition of flow according to the vertical subspaces (involutive) and the horizontal subspace established by this connection. The solution flow (global in $G$ up to lifetime of $Z_t$) is given by left action:
\[
\varphi_t (x)= g_t x,
\]
where $g_t= \exp\{ A\, Z_t \}$. In this example the distributions $\Delta^1$ and $\Delta^2$ in the tangent space $TG$ are given by the horizontal subspaces with respect to the connection $\omega$ and the tangent to the fibres $gH$ (involutive). In order to decompose the flow $\varphi_t$ as in Theorem \ref{theorem: main decomposition}, one has to identify the vector fields $V$ and $h$ as in equation (\ref{Eq: h and V}) in the proof of the theorem, i.e.:
\[
Ax := h + V(\eta, x)
\]
Elements $\eta \in \mathrm{Diff}(\Delta^1, G)$ can be written pointwise (with respect to $x\in G$) as a left action of elements of $G$ at $x$. This action preserves the vertical component, i.e. $g_* \Delta^2 = \Delta^2$ for all $g\in G$. Hence, vector field $V$ above is independent of $\eta$ and one can easily  calculate:
\[
V(x)= \omega (Ax)^*  \ \ \ \     \mbox{ and }  \ \ \ \ h=Ag- \omega (Ax)^*.
\]
By equations (\ref{equation: horizontal components}) and (\ref{Eq: vertical component}) we have that each component of the decomposition $\varphi_t(\cdot)= \eta_t \circ \psi_t (\cdot) $ are given by:
\begin{equation} \label{Eq: eta para fibrado}
d\, \eta_t = R_{\eta_t*} (A\eta_t(\cdot) - \omega (A\eta_t (\cdot))^*)
\end{equation}
and 
\begin{equation} \label{Eq: psi para fibrado}
d\, \psi_t = \mathrm{Ad}(\eta_t) \ \omega (A \eta_t (\cdot))^*.
\end{equation}

Let denote by $g^{H,x}_t \in G$ the $\alpha$-H\"older curve in $G$ such that $g^{H,x}_t x $  is the horizontal lift of $\pi(g_t\, x)$  starting at $x$, i.e. $g^{H,x}_t x $ is horizontal and $g^{H,x}_t x = g_t\, x\, v_t$ for some $v_t \in H$. With this notation, fixing the action at a point $x\in G$, the equations above reduce to well known finite dimensional equations (in G). This is the content of the following

\begin{proposition}
\label{Prop: decomposition in fibre bundles}
Consider the decomposition  $\varphi_t (\cdot)=  \eta_t \circ \psi_t (\cdot)$ of the solution flow of equation (\ref{equation: principal fibre bundle}) according to horizontal and vertical distribution of the fibre bundle in the sense of Theorem \ref{theorem: main decomposition}. Then, at each point $x\in G$, the first component can be written as the left action: 
\[
\eta_t(x)= g^{H,x}_t x, 
\]
and the second component can be written as the right action:
\[
\psi_t(x)= x \, h_t
\]
where $h_t=
x^{-1}\,(g^{H,x}_t)^{-1} \, g_t \, x.
$

\end{proposition}
\proof

The proof of the first equation follows straightforward when one applies equations (\ref{Eq: eta para fibrado})  at a fixed inicial condition $x\in G$: it is the horizontal lift of $\pi (g_t x)$, cf. Theorem  \ref{Thm: Horizontal lift}, using It\^o formula \ref{theorem: Ito formula}. Regarding the second equation of the statement, one sees that $(x^{-1}\,(g^{H,x}_t)^{-1} \, g_t \, x )\in H$ by definition of the horizontal lift: $g^{H,x}_t x = g_t\, x\, v_t$ for some $v_t \in H $. One checks that it solves (\ref{Eq: psi para fibrado}) at a fixed point $x$.

\eop

\bigskip

\noindent \textbf{Trivial fibre bundles:} As a particular case, consider  a trivial principal fibre bundle $\pi : G\times H \rightarrow H$ with strutural group $H$, where $G$ and $H$ are connected Lie groups.  The trivial connection is given by $\omega_{(x,y)}(g_t', h_t')= y^{-1} h_t'\in \mathfrak{h}$. Consider a right invariant YDE in $G \times H$:
\[
d (x_t, y_t) = (A\times B) \  (x_t, y_t)\, dZ_t 
\]
where $A\in \mathfrak{g}$ and $B \in \mathfrak{h}$, the Lie algebras of $G$ and $H$ respectively, with an initial condition $(x_0, y_0)$. Since the connection in this case is invariant by left action of $G\times \{1d\}$, the factor $g^{H,x}_t\in G \times H$ of Proposition 
\ref{Prop: decomposition in fibre bundles} does not depend on $(x,y)$. One recovers the trivial components of the decomposition. In fact we get a global decomposition where the first component is given by the left action:
\[
\eta_t(\cdot, \cdot)= (\exp (AZ_t), 1d) (\cdot,\cdot). 
\]
And the second (vertical) component is given in terms of the right action:
\[
\psi_t(\cdot ,\cdot)= (\cdot , \cdot ) (1d, h_t)
\]
where $h_t=
y^{-1}\, \exp (BZ_t) \, y, 
$ according to Proposition \ref{Prop: decomposition in fibre bundles}.


\bigskip


\begin{thebibliography}{}


\bibitem{Abraham} Abraham, R. ; Marsden, J. E. and Ratiu, T. \textit{Manifolds, tensor analysis, and applications.} Second edition. Applied Mathematical Sciences, 75. Springer-Verlag, New York, 1988.

\bibitem{Applebaum-Estrade} Applebaum, D. and Estrade, A -- Isotropic  L\'evy processes on Riemaniann manifolds. Ann. Probab. \textbf{28}(1),pp. 166-184, 2000.

\bibitem{Bailleul} Bailleul, I. --{\it On the definition of a solution to a rough differential equation.}
 Annales de la Facult{\'e} des Sciences de Toulouse. Math{\'e}matiques, 2019.

\bibitem{Bismut} Bismut, J.-M. -- Mècanique aléatoire, {\it Lectures notes in Mathematics}, 866.
{\it Springer - Verlag}, Berlin - New York, 1981.

\bibitem{Lyons1} Lyons, T. ; Caruana, M. and Levy, T. --
Differential Equations Driven by Rough Paths. Ecole d'Et\'e de Probabilit\'es
de Saint-Flour XXXIV, {\it Lecture Notes in Mathematics} 1908.  Springer, $2007$.

\bibitem{RP}  Castrequini, R. and Catuogno, P.-- A generalized change of variable formula for the Young integral. To appear in {\it Chaos, Solitons and Fractals}, 2022.

\bibitem{CR}  Castrequini, R. and Russo, F. -- Path dependent equations driven by H{\"o}lder processes. {\it Stochastic Analysis and Applications}, v. 37, 3, pp. 480--498, 2019.

\bibitem{Catuogno}  Catuogno, P., Silva, F. and Ruffino, P. -- Decomposition of stochastic flows in manifolds with complementary distributions. {\it Stoch. Dyn.}, v. 13, 4, p. 1350009, 2013.


\bibitem{Colonius and Ruffino} Colonius, F.  and Ruffino, P. --
Nonlinear Iwasawa decomposition of control flows.
{\it Discrete Contin. Dyn. Syst.} v. 18,  2-3, 339--354, 2007.

\bibitem{Cong} Cong, N.; Duc, L. and Hong, P. -- Nonautonomous Young differential equations revisited. {\it J. Dynam. Differential Equations}, v. 30, 4, pp. 1921--1943, 2018.


\bibitem{Manfredo}  do Carmo, M.-- Riemannian geometry. Translated from the second Portuguese edition by Francis Flaherty.  Birkh\"auser Boston, Inc., Boston, MA, 1992.

\bibitem{Elworthy} Elworthy, D. -- Geometric Aspects of Diffusions on Manifolds. Springer, Berlin, 1988.



\bibitem{Friz-Hairer} Friz, P. and Hairer, M. -- {\it A course on rough paths.} 2nd. Edition.
Springer, 2020.

\bibitem{Gubinelli} Gubinelli, M.; Lejay, A. and Tindel, S. -- Young integrals and SPDEs. {\it J. Potential Analysis}, v. 25, 4, pp. 307--326, 2006.


\bibitem{Ikeda-Manabe} Ikeda, N. and Manabe, S. -- Integral of Differential Forms along the Path of Diffusion Processes. {\it Research Institute for Mathematical Sciences}, v.15, pp. 827--852, 1979.




\bibitem{HK} Kunita, H. -- Stochastic Flows and Stochastic Differential Equations - {\it Cambridge University
Press, Cambridge}, 1997.

\bibitem{kobayashi}  Kobayashi, S. and Nomizu, K. -- {\it Foundations of differential geometry}. Vol 1,{\it Wiley} New York, London, 1963.


\bibitem{Ledesma}  Ledesma, D. and  Borges da Silva, F. -- Decomposition of stochastic flow and an averaging principle for slow perturbations. {\it Dyn. Syst.} v. 35,  4, pp. 625--654, 2020.


\bibitem{Lejay} \textrm{Lejay, A.} --Controlled differential equations as Young integrals: a simple approach.
{\it J. Differential Equations}  v. 249,  8, pp. 1777--1798, 2010. 

\bibitem{Lyons}  Li, X. and Lyons, T. --Smoothness of Itô maps and diffusion processes on path spaces. I. {\it
Ann. Sci. École Norm. Sup.} v. 39, 4, 649--677, 2006.





\bibitem{Liao} Liao, M. -- Decomposition of stochastic flows and Lyapunov exponents. {\it Probab.
Theory Rel. Fields} v.117 pp. 589--607, 2000 .


\bibitem{Liao1} Liao, M. -- Invariant Markov processes under Lie group actions. {\it Springer, Cham.}, 2018. 


\bibitem{Malliavin}  Malliavin, M. and  Malliavin, P. -- Factorisations et lois limites de la diffusion horizontale au-dessus d'un espace riemannien sym\'etrique.  Th\'eorie du potentiel et analyse harmonique.  Lecture Notes in Mathematics, v. 404,  pp. 164--217, Springer, Berlin, 1974. 

\bibitem{Melo1} Melo, A., Morgado, L. and Ruffino, P. -- Topology of foliations and decomposition of stochastic flows of diffeomorphisms. {\it Journal of Dynamics and Differential Equations}, v. 30, pp.39--54, 2018.

\bibitem{Melo2} Melo, A., Morgado, L. and Ruffino, P. -- Decomposition of stochastic flows generated by Stratonovich SDEs with jumps. {\it Discrete Contin. Dyn. Syst. - B}, v. 21, 9. pp.3209-–3218, 2016.

\bibitem{Neeb} Neeb, K.-H. -- {\it Infinite Dimesional Lie Groups}. Monastir Summer Schoool, 2009. Available from
hal.archives-ouvertes.fr/docs/00/39/.../CoursKarl-HermannNeeb.pdf.

\bibitem{Omori}  Omori, H. -- Infinite-dimensional Lie groups. v.158, {\it American Mathematical Soc.}, 2017.


\bibitem{Ruzmaikina} Ruzmaikina, A. -- Stieltjes integrals of H{\"o}lder continuous functions with applications to fractional Brownian motion. {\it J. Statist. Phys.}, v. 100, 5, pp. 1049--1069, 2000.


\bibitem{shigekawa}  Shigekawa, I. -- On stochastic horizontal lifts. {\it Z. Wahrsch. Verw. Gebiete}, v.59, pp. 211--221, 1982.

\bibitem{Young} Young, L. -- An inequality of Holder type connected with Stieljes integration. {\it Acta Math}, v.67, pp. 251--258, 1936.





\end{thebibliography}
\end{document}